\documentclass[10pt]{amsart}

\usepackage{amssymb}
\usepackage{amsmath}

\usepackage[colorlinks=true, pdfauthor={David Fithian},
pdfsubject={Congruence Obstructions to Pseudomodularity of Fricke Groups},
pdfkeywords={cusps, Fuchsian groups}]{hyperref}

\newtheorem{theorem}{Theorem}[section]
\newtheorem{proposition}[theorem]{Proposition}
\newtheorem{corollary}[theorem]{Corollary}

\theoremstyle{definition}
\newtheorem{definition}[theorem]{Definition}

\newtheorem{question}{Question}

\newcommand{\PSL}{\operatorname{PSL}}
\newcommand{\PGL}{\operatorname{PGL}}
\newcommand{\fieldQ}{\mathbb{Q}}
\newcommand{\fieldR}{\mathbb{R}}
\newcommand{\ringZ}{\mathbb{Z}}
\newcommand{\ringZhat}{\widehat{\mathbb{Z}}}
\newcommand{\prj}{\mathbb{P}}
\newcommand{\adeles}{\mathbb{A}_{\fieldQ,f}}
\newcommand{\cusps}{\mathcal{C}}

\begin{document}

\title{Congruence obstructions to pseudomodularity of Fricke groups}
\author{David Fithian}
\address{Department of Mathematics\newline
  \indent University of Pennsylvania\newline
  \indent Philadelphia, PA 19104, USA}
\email{fithian@math.upenn.edu}

%\keywords{}
%\date{19 July 2007}

\dedicatory{}

\begin{abstract}
  A pseudomodular group is a finite coarea nonarithmetic Fuchsian group whose cusp set is
  exactly $\prj^1(\fieldQ)$. Long and Reid constructed finitely many of these by
  considering Fricke groups, i.e., those that uniformize one-cusped tori. We prove that a
  zonal Fricke group with rational cusps is pseudomodular if and only if its cusp set is
  dense in the finite adeles of $\fieldQ$. We then deduce that infinitely many such Fricke
  groups are not pseudomodular.
\end{abstract}

\maketitle

\section{Introduction}

A \emph{cusp} of a Fuchsian group $\Gamma\subset\PSL_2(\fieldR)$ is an
$x\in\prj^1(\fieldR)$ that is the unique fixed point of an element of $\Gamma$
(see~\cite{Sh} Ch. 1). The modular group $\PSL_2(\ringZ)$ is a finite coarea Fuchsian
group whose cusps are $\prj^1(\fieldQ) = \fieldQ \cup \{\infty\}$. In~\cite{LR}, Long and
Reid show that there exist finite coarea Fuchsian subgroups of $\PSL_2(\fieldQ)$ that are
\emph{not} commensurable with $\PSL_2(\ringZ)$ (i.e., not arithmetic) and whose cusp set
equals $\prj^1(\fieldQ)$. They call such groups \emph{pseudomodular}.

Long and Reid studied a particular family $\Delta(u^2,2t)$ of Fricke groups as candidates
for pseudomodularity. Fricke groups uniformize one-cusped hyperbolic tori and all such
tori are uniformized by Fricke groups; see~\cite{Ab} and~\cite{Ke}. Among Long and Reid's
stated open problems in~\cite{LR} is the determination of the values
$(u^2,2t)\in\fieldQ\times\fieldQ$ for which $\Delta(u^2,2t)$ is pseudomodular. Recall
that if $\ringZhat$ is the profinite completion of $\ringZ$, then $\adeles =
\fieldQ\otimes_\ringZ\ringZhat$ is an additive topological group having a basis of open
neighborhoods of 0 consisting of $m\ringZhat$ for $m\in\fieldQ$. Then our chief result
is:

\begin{theorem}\label{thm:ad}
  The Fuchsian group $\Delta(u^2,2t)$ is pseudomodular or arithmetic if and only if its
  finite cusps are dense in the ring $\adeles$ of finite adeles over $\fieldQ$.
\end{theorem}

\noindent Theorem~\ref{thm:ad} is in fact true for arbitrary zonal Fuchsian subgroups of
$\PSL_2(\fieldQ)$. We remark that for a given $t$ only finitely many $u^2$ yield
arithmetic groups. We shall use refinements of Theorem~\ref{thm:ad} to give explicit,
infinite families of $\Delta(u^2,2t)$ whose cusps are proper subsets of
$\prj^1(\fieldQ)$. For example:

\begin{theorem}\label{thm:infnonpsm}
  Let $p$ be a prime and $t$ an integer at least 2. Then $\Delta(p^{-2},2t)$ is neither
  pseudomodular nor arithmetic. In particular, infinitely many $\Delta(u^2,2t)$ are
  neither pseudomodular nor arithmetic.
\end{theorem}

Our method is to find $\Delta(u^2,2t)$-invariant subsets of $\fieldQ$ that are defined
number-theoretically, either adelically or $p$-adically. In~\cite{LR}, Long and Reid
exhibit finitely many $\Delta(u^2,2t)$ that are neither pseudomodular nor arithmetic. For
each such group, they provide a rational number fixed by a hyperbolic element of
$\Delta(u^2,2t)$. Such fixed points cannot also be cusps; see~\cite{Be}, p. 199. We do not
know whether rational hyperbolic fixed points exist for all non-pseudomodular $\Delta$,
and in any case, our proofs do not require or produce them. So far, our results disprove
pseudomodularity with a finite amount of data, in terms of either finitely many primes or
(finite-index) congruence subgroups of the modular group $\PSL_2(\ringZ)$. In the final
section of this paper we discuss some further questions and possible generalizations.

\section{Results}

The numerator and denominator of a rational number will be taken coprime with the
denominator positive. Let $v_p:\fieldQ^\times\to\ringZ$ be the usual discrete valuation at
prime $p$. If $p$ does not divide the denominator of a rational number $x$ (i.e.,
$v_p(x)\geq 0$) then say $x$ is \emph{integral} at $p$. Integral $x$ have representatives
in the ring of $p$-adic integers $\ringZ_p$ and we write $x\equiv m\ (\bmod\ p)$ to mean
that the residue of $x$ mod $p\ringZ_p$ is $m$ in $\ringZ_p/p\ringZ_p\cong
\ringZ/p\ringZ$. Finally, for nonzero rational $x$, we say $p\mid x$ or $p$ \emph{divides}
$x$ when $v_p(x) > 0$.

Denote by $\cusps_\infty(G)$ the cusp set of a Fuchsian group $G$ and let the finite cusps
be denoted by $\cusps(G):=\cusps_\infty(G)\setminus\{\infty\}$. For rationals $u^2$ and
$t$ with $0 < u^2 < t-1$, let $\Delta(u^2,2t)$ be the subgroup of $\PSL_2(\fieldR)$
generated by the hyperbolic elements
\[
g_1 = \frac{1}{\sqrt{-1+t-u^2}}\left(
  \begin{array}{cc}
    t-1 & u^2 \\
    1 & 1 \\
  \end{array}
\right),\quad g_2 = \frac{1}{u\sqrt{-1+t-u^2}}\left(
  \begin{array}{cc}
    u^2 & u^2 \\
    1 & t-u^2 \\
  \end{array}
\right).
\]
As in~\cite{LR}, $\Delta$ is a zonal Fricke group freely generated by $g_1$ and $g_2$
with $\cusps_\infty(\Delta)$ a nonempty subset of $\prj^1(\fieldQ)$. In fact, by
considering the traces of $g_1$, $g_2$ and $g_1g_2$ and using results in \S1 of~\cite{CS}
on varieties of group representations, we can show that every Fricke group with cusps in
$\prj^1(\fieldQ)$ is conjugate in $\PGL_2(\fieldQ)$ to some $\Delta(u^2,2t)$. Therefore,
we can study all such Fricke groups by looking at the groups $\Delta(u^2,2t)$. We are
interested in when $\cusps_\infty(\Delta)=\prj^1(\fieldQ)$ or, equivalently, when
$\cusps(\Delta)=\fieldQ$. Let $\Lambda(u^2,2t)$ be the kernel of the homomorphism
$\Delta(u^2,2t) \to \ringZ/2\oplus\ringZ/2$ given by $g_1\mapsto (1,0)$ and $g_2\mapsto
(0,1)$. Then $\Lambda(u^2,2t)\subset\PSL_2(\fieldQ)$ and $\cusps(\Lambda) =
\cusps(\Delta)$. If $\cusps(\Lambda)$ is not dense in $\adeles$, or in some finite
product $\prod_i\fieldQ_{p_i}$ with $\fieldQ$ embedded diagonally, then
$\cusps(\Lambda)\neq\fieldQ$ since $\fieldQ$ is dense in each of those sets.

Given this obstruction to pseudomodularity, we ask three questions: (i) when are the
finite cusps of $\Delta(u^2,2t)$ dense in a specified finite product
$\prod_i\fieldQ_{p_i}$, (ii) is cusp density in all such products a sufficient condition
for pseudomodularity, and (iii) how do the answers to the former two questions vary as
$u^2$ and $t$ are themselves varied $p$-adically? The propositions below
address these questions and, in particular, prove Theorem~\ref{thm:infnonpsm}.

\begin{proposition}\label{prop:square}
  Let $p$ be prime. If $v_p(t)\geq 0$ and $v_p(u^2)\leq -2$, or if $v_p(t) < 0$ and
  $v_p(u^2)\leq 2(v_p(t)-1)$, then $\cusps(\Delta(u^2,2t))$ is not dense in $\fieldQ_p$.
\end{proposition}

\noindent Considering conditions at two primes, we can show:

\begin{proposition}\label{prop:2primes}
  If $p$ and $q$ are prime, $v_p(u^2) = -1 = v_q(u^2)$, and $t$ is integral at $p$ and at
  $q$, then $\cusps(\Delta(u^2,2t))$ is not dense in $\fieldQ_p\times\fieldQ_q$.
\end{proposition}

\noindent Results similar to Proposition~\ref{prop:2primes} hold for $t$ that are
non-integral at $p$ or at $q$; we omit their statements for brevity. As a corollary to
the above two propositions, whenever $t$ is an integer and the denominator of $u^2$ is
composite, $\Delta(u^2,2t)$ is not pseudomodular. We have an even stronger result for
integral $t$:

\begin{proposition}\label{prop:tisint}
  Let $t$ be an integer and suppose $\Delta(u^2,2t)$ is pseudomodular. Then
  \begin{enumerate}
  \item[\emph{(a)}] $u^2$ has prime or unit denominator, say $p$,
  \item[\emph{(b)}] if this $p$ is an odd prime, then $p$ does not divide $t$, and
  \item[\emph{(c)}] for all odd primes $q$ dividing $t$, $u^2$ (necessarily in $\ringZ_q$)
    is equivalent to $0$ or $-1\bmod q$.
  \end{enumerate}
\end{proposition}

To prove each of the above propositions, we give a proper nonempty $\Delta$-invariant
subset $U$ of $\fieldQ$ that is open in the topology induced by that of $\fieldQ_p$ or
$\fieldQ_p\times\fieldQ_q$. For example, under the hypotheses of
Proposition~\ref{prop:2primes}, the set of nonzero rationals $x$ for which exactly one of
$v_p(x)$ and $v_q(x)$ is negative is $\Delta$-invariant. To prove that our $U$ are
$\Delta$-invariant, we show that the generators $g_i^{\pm 1}$ of $\Delta$ map $U$ into
itself.

The difficulty of our approach lies in finding such sets $U$. To do this, we search
computed data for patterns in valuations and congruences. We build a data set by taking a
tuple $(x_1,\ldots,x_m)$ of rational numbers and applying to it several elements of
$\Delta$. Useful values for the $x_i$ are rational hyperbolic (or \emph{special}) fixed
points, such as those tabulated in~\cite{LR}. In some cases, a tuple of length one yields
some information, but typically we let the tuple be given by multiple pairs of special
fixed points, with each pair fixed by a single hyperbolic element of $\Delta$. When we
have no data on special fixed points, we select the $x_i$ with particular valuation or
congruence properties. In any case, a successful search suggests a candidate $U$, which we
then confirm as described above. Varying parameters such as the denominator of $u^2$, in
such a way as to preserve the proof of the $\Delta$-invariance, allows us to eliminate
infinite families of candidates for pseudomodularity and divorce our results from computed
data.

In the other direction, to show that $\cusps(\Delta)$ is dense in a finite product
$H=\prod_i\fieldQ_{p_i}$, we show that for each rational $x$, we can move $x$ arbitrarily
close to the cusp 0 with elements of (a group commensurable with)
$\Delta$. We use this argument to prove that for $t$ a prime integer, we have identified
above \emph{all} cases for which $\cusps(\Delta(u^2,2t))$ is not dense in some $H$:

\begin{proposition}
  If $t$ is prime, $u^2$ has prime denominator not equal to $t$ and $u^2\equiv 0$ or
  $-1\bmod t$, then $\cusps(\Delta(u^2,2t))$ is dense in every finite product
  $\prod_i\fieldQ_{p_i}$ and hence is dense in the product $\prod_p\fieldQ_p$ over all
  primes.
\end{proposition}

\noindent There are groups with special fixed points to which this proposition applies,
such as $\Delta(6/11,6)$ with a special fixed point of $1/4$. Consequently:

\begin{corollary}\label{cor:density}
  Density of $\cusps(\Delta(u^2,2t))$ in the product $\prod_p\fieldQ_p$ of all $p$-adic
  fields is not a sufficient condition for pseudomodularity.
\end{corollary}

\noindent Both this corollary and Theorem~\ref{thm:ad} can be interpreted in terms of congruence
data with respect to all primes. While density in the adelic topology implies
pseudomodularity, density in the product topology on all $p$-adic fields does not.
Theorem~\ref{thm:ad} holds because an orbit of a rational under a fixed translation
$x\mapsto x + t$ with $t\in\fieldQ$ is open in $\fieldQ$ in the topology induced by
$\adeles$.

\section{Questions}

Since $\Lambda(u^2,2t)$ is a finitely generated subgroup of $\PSL_2(\fieldQ)$, it is
always a subgroup of $\PSL_2(\ringZ_S)$, where $\ringZ_S =
\ringZ[p_1^{-1},\ldots,p_r^{-1}]$ for some minimal set of primes $p_1,\ldots,p_r$. By
Riemann-Hurwitz, $\Lambda(u^2,2t)$ always has exactly four orbits of cusps, so if we can
show that $\Lambda(u^2,2t)$ has more than four orbits in its action on $\prj^1(\fieldQ)$,
then some rationals cannot be cusps, whence $\Lambda(u^2,2t)$ is not pseudomodular.

\begin{question}
  For a non-pseudomodular $\Lambda(u^2,2t)\subset\PSL_2(\ringZ_S)$ as above, can we find a
  subgroup $K$ of (finite index in) $\PSL_2(\ringZ_S)$ such that $\Lambda(u^2,2t)\subset
  K$ and $K$ has more than four orbits in its action on $\prj^1(\fieldQ)$?
\end{question}

For the groups already eliminated by our earlier propositions, the answer is yes. We also
have a positive answer for some groups whose cusps are dense in the product of all
$\fieldQ_p$'s, such as $\Lambda(6/11,6)$. In fact, we can construct a
$K\subset\PSL_2(\ringZ[2^{-1}])$ containing $\Lambda(6/11,6)$ that has eight orbits in its
action on $\prj^1(\fieldQ)$. An explicit description of the orbits of $K$ gives us a
$\Delta(6/11,6)$-invariant, nonempty, proper subset $X$ of $\prj^1(\fieldQ)$ that is open
in the 33-congruence topology, as defined here.

\begin{definition}\label{def:congtop}
  Denote by $\Gamma(M)$ the principal congruence subgroup of $\PSL_2(\ringZ)$ of level
  $M$.  For a finite set of primes $p_1,\ldots,p_r$, let $N=p_1\ldots p_r$ and let
  $\mathcal{V}_N$ be the set of orbits $\Gamma(N^j)\cdot 0$ as $j$ ranges over positive
  integers. The \emph{$N$-congruence topology} on $\prj^1(\fieldQ)$ is the topology
  generated by the $\PSL_2(\ringZ)$-translates of $\mathcal{V}_N$.
\end{definition}

\noindent If we replace $\Gamma(N^j)$ with $\Gamma^0(N^j)$ (matrices lower-triangular mod
$N^j$) in this definition then we recover the (coarser) topology on $\prj^1(\fieldQ)$
induced by inclusion in $\prod_{i=1}^r \prj^1(\fieldQ_{p_i})$. Thus, in
Propositions~\ref{prop:square} and \ref{prop:2primes}, we also have cusp sets failing to
be dense in some congruence topology.

\begin{question}
  If the cusps of $\Delta$ are dense in $\prj^1(\fieldQ)$ in all congruence topologies, is
  $\Delta$ necessarily pseudomodular or arithmetic? Equivalently, if $\Delta$ is neither
  pseudomodular nor arithmetic, is there a congruence subgroup $\Gamma$ of
  $\PSL_2(\ringZ)$ so that $\cusps(\Delta)$ misses some orbit $\Gamma\cdot
  x\in\prj^1(\fieldQ)$?
\end{question}

More generally, we can allow non-congruence subgroups $\Gamma$ and we can ask analogous
questions for larger families of Fuchsian groups or for Kleinian groups whose cusps lie in
an imaginary quadratic number field.

\end{document}